%&amstex

\input amsppt.sty

\def\PSH{\operatorname{PSH}}
\def\phi{\varphi}
\def\SH{\operatorname{SH}}
\def\Re{\operatorname{Re}}
\NoBlackBoxes

\document
\topmatter
\title
Effective formulas for invariant functions --
case of elementary Reinhardt domains
\endtitle
\rightheadtext{Effective formulas for invariant functions}
\author
Peter Pflug \& W\l odzimierz Zwonek
\endauthor
\address
Carl von Ossietzky Universit\"at Oldenburg, Fachbereich Mathematik, Postfach 2503,
D-26111 Oldenburg, Germany
\endaddress
\email
pflug\@mathematik.uni-oldenburg.de
\endemail
\address
Uniwersytet Jagiello\'nski, Instytut Matematyki, Reymonta 4,
30-059 Krak\'ow, Poland
\endaddress
\email zwonek\@im.uj.ed.pl
\endemail
\abstract
In the paper we find effective formulas for the invariant functions,
appearing in the theory of several complex variables,
of the elementary Reinhardt domains. This gives us the first example
of a large family of domains for which the functions are calculated
explicitly.
\endabstract
\thanks The research started while the second author was visiting
Universit\"at Oldenburg, whose stay was enabled by Volkswagen
Stiftung Az. I/71 062.
The second author was also supported by KBN Grant No 2 PO3A 060 08
\endthanks
\endtopmatter

\subheading{0. Introduction}  Holomorphically invariant functions  and
pseudometrics have  proved to  be very  useful in  the theory  of several  complex variables.
Nevertheless, the  problem of finding effective formulas for
the objects has turned out to be very difficult.
So far  there have been very few examples of domains
for  which the  formulas for  these  functions  are known  explicitly.

Among many different invariant functions and pseudometrics
let us mention the Lempert and Green functions,
the Kobayashi and Carath\'eodory pseudodistances
as well as their infinitesimal versions i.e.
the Kobayashi--Royden, Carath\'eodory and Azukawa pseudometrics.

Due  to
Lempert's  theorem (see \cite{L 1,2}) all  holomorphically invariant  functions and  pseudometrics
coincide in  the class of  convex domains, therefore, these are
the non--convex domains, which
may deliver us a great deal of different invariant functions, not the
convex ones.
But even in convex case
it is  difficult to find  explicit formulas for the objects involved.
Among few results in this direction let us
mention here the special case of convex (see \cite{BFKKMP},
\cite{JP 2}) and non--convex (see \cite{PZ})
ellipsoids.
The other class of non--convex domains for which some of the functions were
calculated  is a class  of elementary Reinhardt
domains (see \cite{JP 1,2}).
In our paper we extend the results obtained in those domains
for all invariant functions and pseudometrics mentioned earlier.
The formulas obtained enable us to understand better the mutual relations
between the invariant objects and give surprising solutions to some problems.

\subheading{1. Definitions, notations and main results} By $E$ we will
always denote the unit disc in $\Bbb C$. We put
$m(\lambda_1,\lambda_2):=\frac{|\lambda_1-\lambda_2|}
{|1-\bar\lambda_1\lambda_2|}$
for $\lambda_1,\lambda_2\in E$ and
$\gamma_E(\lambda;\alpha):=\frac{|\alpha|}{1-|\lambda|^2}$,
$\lambda\in   E$, $\alpha\in\Bbb C$.

Let $D$ be a domain in
$\Bbb C^n$. Following \cite{L 1}, \cite{Ko 1}, \cite{Kl 1,2},
\cite{C}, \cite{A}, \cite{R} and \cite{JP2} for $(w,z)\in D$ and
$(w;X)\in D\times \Bbb C^n$ we define the following functions:
$$
\align
\tilde k_D^*(w,z)&:=
\inf\{m(\lambda_1,\lambda_2):\exists\phi\in\Cal    O(E,D),\;
\phi(\lambda_1)=w,\; \phi(\lambda_2)=z\},\\
k_D^*(w,z)&:=\tanh k_D(w,z),\\
\text{ where $k_D$ is the largest } &
\text{ pseudodistance smaller or equal than $\tilde k_D:=
\tanh^{-1}\widetilde k^*_D$},\\
g_D(w,z)&:=
\sup\{u(z):\\
\log u\in\PSH(D,[-\infty,0)),&\;\exists M,R>0:
u(\zeta)\leq M||\zeta-w||\text{ for } \zeta\in D,\;
||\zeta-w||<R\},\\
c_D^*(w,z)&:=\sup\{m(\phi(w),\phi(z)):\;\phi\in\Cal O(D,E)\};
\endalign
$$
and also their infinitesimal versions:
$$
\align
\kappa_D(w;X)&:=
\inf\{\gamma_E(\lambda;\alpha):\exists \phi\in\Cal O(E,D),\;
\phi (\lambda)=z,\alpha
\phi^{\prime}(\lambda)=X\},\\
A_D(w;X)&:=
\limsup_{\lambda\not\to 0}\frac{g_D(w,w+\lambda X)}{|\lambda|},\\
\gamma_D(w;X)&:=
\sup\{\gamma_E(\phi(w);\phi^{\prime}(w)X):\;\phi\in\Cal O(D,E)\}.
\endalign
$$
The function $\tilde k_D^*$ (respectively, $g_D$, $k_D^*$, $c_D^*$) is called
{\it the Lempert function} (respectively, {\it the Green function},
{\it  the   Kobayashi }  and  {\it Carath\'eodory   pseudodistance}).
The  function
$\kappa_D$  (respectively,   $A_D$  and  $\gamma_D$) is called  {\it   the
Kobayashi--Royden}
(respectively, {\it Azukawa } and {\it Carath\'eodory--Reiffen})
pseudometric.

Note that the functions $\tilde k_D^*$, $k_D^*$ and $c_D^*$
are always symmetric, whereas $g_D$ need not have the property.
For the  basic properties of the functions  defined we refer the
interested  reader to  \cite{JP 2}.
Let  us mention  here only  some basic
relations between the objects involved:
$$
\gather
\tilde k_D^*\geq k_D^*\geq c_D^*,\quad \tilde k_D^*\geq g_D\geq c_D^*,\\
\kappa_D\geq A_D\geq \gamma_D.
\endgather
$$
A mapping $\phi\in\Cal O(E,D)$ is called a \it $\tilde k_D$-geodesic \rm
for $(w,z)$, $w\neq z$ if $\phi(\lambda_1)=w$, $\phi(\lambda_2)=z$ and
$m(\lambda_1,\lambda_2)=\tilde k_D^*(w,z)$ for suitable
$\lambda_1,\lambda_2\in E$.

The class of domains we are intersted in is defined below.

For $\alpha=(\alpha_1,\ldots,\alpha_n)\in\Bbb R_+^n$, $n>1$,
($\Bbb R_+:=(0,\infty)$) define
$$
D_{\alpha}:=\{z\in\Bbb C^n:|z_1|^{\alpha_1}\cdot\ldots\cdot|z_n|^{\alpha_n}
<1\}.
$$
We say that $\alpha$ is of a \it rational type \rm if
there are $t>0$, $\beta=(\beta_1,\ldots,\beta_n)\in \Bbb N_*^n$ such
that $\alpha=t\beta$. We say that $\alpha$ is of {\it irrational type }
if $\alpha$  is not of  a rational type.  Remark that in case when
$\alpha$ is  of
rational type we  may without loss of generality  assume that all
$\alpha_j$'s are relatively prime natural numbers.
We define also
$$
\tilde D_{\alpha}:=\{z\in D_{\alpha}:z_1\cdot\ldots\cdot z_n\neq 0\},
$$
If $\alpha\in\Bbb N_*^n$, then we denote
$$
\gather
z^{\alpha}:=z_1^{\alpha_1}\cdot\ldots\cdot z_n^{\alpha_n},\quad
F^{\alpha}(z):=z^{\alpha},\\
F^{\alpha}_{(r)}(z)X:=\sum_{\beta_1+\ldots+\beta_n=r}
\frac{1}{\beta_1!\cdot
\ldots\cdot\beta_n!}\frac{\partial^{\beta_1+\ldots+\beta_n}F^{\alpha}(z)}
{\partial z_1^{\beta_1}\ldots \partial z_n^{\beta_n}}X^{\beta},
\quad z,X\in\Bbb C^n.
\endgather
$$
Note that the domain $D_{\alpha}$  is always unbounded, Reinhardt, complete,
and pseudoconvex but not convex.

As  mentioned in  Introduction some   of the  invariant functions  for domains
$D_{\alpha}$ are explicitly known.
We gather  the results  known so  far in the following
theorem
\proclaim{Theorem  1 {\rm (see \cite{JP 2})}}
If  $\alpha\in\Bbb  N_*^n$,   where  $\alpha_j$'s  are
relatively prime, then :
$$
\align
c_{D_{\alpha}}^*(w,z)&=m(w^{\alpha},z^{\alpha}),\\
g_{D_{\alpha}}(w,z)&=m(w^{\alpha},z^{\alpha})^{\frac{1}{r}},\\
\gamma_{D_{\alpha}}(w;X)&=
\gamma_E(w^{\alpha},(F^{\alpha})^{\prime}(w)X),\\
A_{D_{\alpha}}(w;X)&=\left(\gamma_E(w^{\alpha},F^{\alpha}_{(r)}(w)X)
\right)^{\frac{1}{r}},\qquad (w,z)\in D_{\alpha}\times D_{\alpha},
(w;X)\in D_{\alpha}\times \Bbb C^n,
\endalign
$$
where $r$ is the order of vanishing of the function
$F^{\alpha}(\cdot)-F^{\alpha}(w)$ at $w$.

If $\alpha$ is of irrational type, then
$$
\align
c_{D_{\alpha}}^*(w;X)&=0,\\
\gamma_{D_{\alpha}}(w;X)&=0,\qquad
(w,z)\in D\times D, (w;X)\in D\times \Bbb C^n.
\endalign
$$
\endproclaim
In our paper we extend the results of Theorem 1
to other invariant functions and
pseudometrics and  we find the  remaining formulas for  the Green function
(and Azukawa pseudometric) in
the irrational  case. The  results are  presented  in  two theorems.
One of  them
concerns with rational, while the other  one with irrational $\alpha$.
In both
theorems in case  of the Lempert function
the formulas may seem  to be incomplete
(not all the cases are
covered);  nevertheless, because of the symmetry of  both functions
one easily obtains the formulas in remaining cases.

\proclaim{Theorem 2} Assume that $\alpha\in\Bbb N_*^n$ with $\alpha_j$'s relatively
prime.
Let $(w,z)\in D_{\alpha}\times D_{\alpha}$, $(w;X)\in
D_{\alpha}\times\Bbb C^n$.
Denote $\Cal J:=\{j\in\{1,\ldots,n\}:w_j=0\}=\{j_1,\ldots,j_k\}$.
Then we have
$$
\gather
\tilde k_{D_{\alpha}}^*(w,z)=\\
\cases
\min\{m(\lambda_1,\lambda_2):\lambda_1,
\lambda_2\in E,
\lambda_1^{\min\{\alpha_j\}}=w^{\alpha},
\lambda_2^{\min\{\alpha_j\}}=z^{\alpha}\},&\text{        if       $w,z\in\tilde
D_{\alpha}$}\\
|z^{\alpha}|^{\frac{1}{\alpha_{j_1}+\ldots+
\alpha_{j_k}}},&\text{ if $\Cal J\neq\emptyset$}
\endcases;
\endgather
$$
$$
k_{D_{\alpha}}^*(w,z)=\min\{m((w^{\alpha})^{\frac{1}{\min\{\alpha_j\}}},
(z^{\alpha})^{\frac{1}{\min\{\alpha_j\}}})\},
$$
where the minimum is taken over all possible roots;

in the infinitesimal case we have
$$
\multline
\kappa_{D_{\alpha}}(w;X)=\\
\cases
\gamma_E((w^{\alpha})
^{\frac{1}{\min\{\alpha_k\}}},(w^{\alpha})
^{\frac{1}{\min\{\alpha_k\}}}\frac{1}{\min\{\alpha_k\}}
\sum_{j=1}^n\frac{\alpha_j X_j}
{w_j}),&\text{ if $\Cal J=\emptyset$,}\\
(|w_1|^{\alpha_1}\cdot\ldots\cdot|X_{j_1}|
^{\alpha_{j_1}}\cdot\ldots\cdot|X_{j_k}|^{\alpha_{j_k}}
\cdot\ldots\cdot |w_n|^{\alpha_n})^{\frac{1}
{\alpha_{j_1}+\ldots+\alpha_{j_k}}},&
\text{ if $\Cal J\neq\emptyset$}
\endcases.
\endmultline
$$
\endproclaim
Observe that if $\min\{\alpha_j\}=1$, then $\tilde k_{D_{\alpha}}^*
(w,z)=g_{D_{\alpha}}(w,z)$
for $w,z\in\tilde D_{\alpha}$; otherwise, if $w^{\alpha}\neq z^{\alpha}$,
then the Green function is strictly less than the Lempert function.

\vskip2ex
In the irrational case unlike in the rational one, these are not only
the Lempert function, Kobayashi pseudodistance and Kobayashi--Royden
pseudometric, which have not been calculated so far
but also the Green function and the Azukawa pseudometric.

\proclaim{Theorem 3} Assume that $\alpha$ is of irrational type.
Let $(w,z)\in D_{\alpha}\times D_{\alpha}$, $(w;X)\in
D_{\alpha}\times\Bbb C^n$.
Denote $\Cal J:=\{j\in\{1,\ldots,n\}:w_j=0\}=\{j_1,\ldots,j_k\}$.
Then we have
$$
\gather
\tilde k_{D_{\alpha}}^*(w,z)=\\
\cases
m((|w_1|^{\alpha_1}\cdot\ldots\cdot|w_n|^{\alpha_n})
^{\frac{1}{\min\{\alpha_j\}}},
(|z_1|^{\alpha_1}\cdot\ldots\cdot|z_n|^{\alpha_n})
^{\frac{1}{\min\{\alpha_j\}}}),&
\text{ if $w,z\in\tilde D_{\alpha}$}\\
(|z_1|^{\alpha_1}\cdot\ldots\cdot|z_n|^{\alpha_n})
^{\frac{1}{\alpha_{j_1}+\ldots+\alpha_{j_k}}},&
\text{ if $\Cal J\neq\emptyset$}
\endcases;
\endgather
$$
$$
\align
k_{D_{\alpha}}^*(w,z)&=
m\left(\left(\prod_{j=1}^n|w_j|^{\alpha_j}\right)
^{\frac{1}{\min\{\alpha_j\}}},
\left(\prod_{j=1}^n|z_j|^{\alpha_j}\right)
^{\frac{1}{\min\{\alpha_j\}}}\right),\\
g_{D_{\alpha}}(w,z)&=
\cases
0,&\text{ if $\Cal J=\emptyset$},\\
(|z_1|^{\alpha_1}\cdot\ldots\cdot|z_n|^{\alpha_n})
^{\frac{1}{\alpha_{j_1}+\ldots+\alpha_{j_k}}},&\text{ if $\Cal J\neq \emptyset$}
\endcases;
\endalign
$$
in the infinitesimal case we have:
$$
\gather
\kappa_{D_{\alpha}}(w;X)=\\
\cases
\gamma_E(\left(\prod_{j=1}^n|w_j|^{\alpha_j}\right)
^{\frac{1}{\min\{\alpha_k\}}},\left(\prod_{j=1}^n|w_j|^{\alpha_j}\right)
^{\frac{1}{\min\{\alpha_k\}}}\frac{1}{\min\{\alpha_k\}}
\sum_{j=1}^n\frac{\alpha_j X_j}
{w_j})&,\text{if $\Cal J=\emptyset$,}\\
(|w_1|^{\alpha_1}\cdot\ldots\cdot|X_{j_1}|
^{\alpha_{j_1}}\cdot\ldots\cdot|X_{j_k}|^{\alpha_{j_k}}
\cdot\ldots\cdot
|w_n|^{\alpha_n})^{\frac{1}{\alpha_{j_1}+\ldots+\alpha_{j_k}}}&,
\text{if $\Cal J\neq\emptyset$}
\endcases;
\endgather
$$
$$
\gather
A_{D_{\alpha}}(w;X)=\\
\cases
0&,\text{ if $\Cal J=\emptyset$},\\
(|w_1|^{\alpha_1}\cdot\ldots\cdot|X_{j_1}|^{\alpha_{j_1}}
\cdot\ldots\cdot|X_{j_k}|^{\alpha_{j_k}}
\cdot\ldots\cdot|w_n|^{\alpha_n})
^{\frac{1}{\alpha_{j_1}+\ldots+\alpha_{j_k}}}&,\text{ if $\Cal J\neq \emptyset$}
\endcases .
\endgather
$$
\endproclaim
Observe that for an arbitrary balanced pseudoconvex domain $D$
we always have that $\tilde k_D^*(0,z)=h_D(z)$, $z\in D$,
where $h_D$ denotes the Minkowski function for $D$. In the above formula
we have that $k_{D_{\alpha}}^*(0,z)<h_{D_{\alpha}}(z)$,
$0\neq z\in D_{\alpha}$. It would be interesting to find the general form
of $k_D^*(0,\cdot)$ in the case when $D$
is an arbitrary balanced pseudoconvex domain.

\subheading{2. Auxiliary results}

For $z\in\Bbb C^n$ put
$$
T_z:=\{(e^{i\theta_1}z_1,\ldots,e^{i\theta_n}z_n):\theta_j\in\Bbb R\}.
$$
Note that $T_z$ is a group with the multiplication defined as follows:
$$
(e^{i\theta_1}z_1,\ldots,e^{i\theta_n}z_n)
\circ (e^{i\tilde\theta_1}z_1,\ldots,e^{i\tilde\theta_n}z_n):=
(e^{i(\theta_1+\tilde\theta_1)}z_1,\ldots,e^{i(\theta_n+\tilde\theta_n)}z_n).
$$
Define
$T_{z,\alpha}$ as the subgroup of $T_z$ generated by the set
$$
\{(e^{i\frac{\alpha_{j_1}}{\alpha_1}2k_1\pi}z_1,\ldots,
e^{i\frac{\alpha_{j_n}}{\alpha_n}2k_n\pi}z_n):\;j_1,\ldots,j_n\in
\{1,\ldots,n\},\;k_1,\ldots,k_n\in \Bbb Z\}
$$
Note that if $\alpha$ is of a rational type, then
$T_{z,\alpha}$ is finite; more precisely, if we assume that
$\alpha\in\Bbb N_*^n$ and $\alpha_j$'s are relatively prime, then
$$
T_{z,\alpha}=\{(\varepsilon_1z_1,\ldots,\varepsilon_nz_n),\text{ where
$\varepsilon_j^{\alpha_j}=1$}\}.
$$
However, if $\alpha$ is of irrational type,
then a well--known theorem of Kronecker (see \cite{HW}, theorem 439)
gives that
$$
\overline{T_{z,\alpha}}=T_z.\tag{1}
$$
For $\mu\in E_*$ we define
$$
\Phi_{\mu}:\Bbb C^{n-1}\owns(\lambda_1,\ldots,\lambda_{n-1})\to
(e^{\alpha_n\lambda_1},\ldots,e^{\alpha_n\lambda_{n-1}},
\mu e^{-\alpha_1\lambda_1}\cdot\ldots\cdot e^{-\alpha_{n-1}\lambda_{n-1}})
\in D_{\alpha}
$$
Put
$$
V_{\mu}:=\Phi_{\mu}(\Bbb C^{n-1}),\; \mu\in E_*,\quad V_0:=\{z\in\Bbb C^{n}:
z_1\cdot\ldots\cdot z_n=0\}.
$$
Note that
$$
\bigcup_{\mu\in E}V_{\mu}=D_{\alpha}.
$$

\subheading{Remark 4} Let $\mu\in E_*$. Assume that $w,z\in V_{\mu}$, and
$X\in\Bbb C^n$ fulflis the equality $\sum_{j=1}^n\frac{\alpha_jX_j}{w_j}=0$.
Then
$$
\gather
\tilde k_{D_{\alpha}}^*(w,z)=0,\\
\kappa_{D_{\alpha}}(w;X)=0.
\endgather
$$
In fact, $w=\Phi_{\mu}(\lambda)$, $z=\Phi_{\mu}(\gamma)$ for some
$\lambda,\gamma\in \Bbb C^{n-1}$, so
$$
\tilde k_{D_{\alpha}}^*(w,z)=\tilde k_{D_{\alpha}}^*(\Phi_{\mu}(\lambda),
\Phi_{\mu}(\gamma))\leq \tilde k_{\Bbb C^{n-1}}^*(\lambda,\gamma)=0.
$$
To see the second equality note that assuming
$\Phi_{\mu}(\lambda)=w$ we have
$$
\Phi_{\mu}^{\prime}(\lambda)(Y)=\left[\alpha_nw_1Y_1,\ldots,
\alpha_nw_{n-1}Y_{n-1},-\sum_{j=1}^{n-1}\alpha_jw_nY_j\right],
\quad Y\in\Bbb C^{n-1}.
$$
One may easily verify that
$$
\Phi_{\mu}^{\prime}(\lambda)(\Bbb C^{n-1})=\left\{X\in\Bbb C^n:
\sum_{j=1}^n\frac{\alpha_jX_j}{w_j}=0\right\}.
$$
Note that
$$
0=\kappa_{\Bbb C^{n-1}}(\lambda;Y)\geq
\kappa_{D_{\alpha}}(\Phi_{\mu}(\lambda),
\Phi_{\mu}^{\prime}(\lambda)Y),\quad Y\in \Bbb C^{n-1}
$$
which finishes the proof.

In the proof of Lemma 5 below we shall replace $E$
in the definition of the Lempert function with
$H:=\{x+iy:1>x>-1\}$.

\proclaim{Lemma 5} Fix $w,z\in D_{\alpha}$. Take any $\tilde z\in T_{z,\alpha}$.
Then for any $\phi\in\Cal O(E,D_{\alpha})$ such that
$\phi(\lambda_1)=w$, $\phi(\lambda_2)=z$, $\lambda_1\neq \lambda_2$
there is $\tilde\phi\in \Cal O(E,D_{\alpha})$ such that
$\tilde \phi(\lambda_1)=w$ and $\tilde \phi(\lambda_2)=\tilde z$.

Consequently,
$$
\tilde k_{D_{\alpha}}^*(w,z)=\tilde k_{D_{\alpha}}^*(w,\tilde z)\text{ for
any $\tilde z\in T_{z,\alpha}$}.
$$
\endproclaim
\demo{Proof} For the proof of the lemma it is enough to take any mapping
$$
\phi\in\Cal O(H,D_{\alpha}),\quad\phi(0)=w,\quad\phi(it)=z,\;t>0.
$$
Define the mapping ($k_n\in\Bbb Z$ is fixed)
$$
\tilde\phi:H\owns\lambda\to(\phi_1(\lambda),\ldots,\phi_{n-2}(\lambda),
e^{-2k_n\pi\frac{\lambda}{t}}\phi_{n-1}(\lambda),
e^{\frac{\alpha_{n-1}2k_n\pi\lambda}{\alpha_nt}}\phi_n(\lambda))\in D_{\alpha}.
$$
We have
$$
\tilde\phi(0)=w,\quad\tilde \phi(it)=
(z_1,\ldots,z_{n-1},e^{i\frac{\alpha_{n-1}}{\alpha_n}2k_n\pi}z_n).
$$
Note that we may replace $\alpha_{n-1}$ above with any other $\alpha_j$ and
$z_n$ with $e^{i\frac{\alpha_j}{\alpha_n}2k_n\pi}z_n$,
and also we may continue the procedure as above with the next components $z_j$
to be varied, which would finish the proof.
\qed
\enddemo

\subheading{Remark 6} From the proof of Lemma 5 we have also the following
property:

Fix $\alpha\in \Bbb N_*^n$, $\alpha_j$'s relatively prime and
$0<\delta_1\leq m(\lambda_1,\lambda_2)\leq\delta_2<1$. Take any
$\psi\in\Cal O(E,\Bbb C^n)$, $\psi(E)\subset\subset (\Bbb C_*)^n$ and
choose $z\in\Bbb C^n$ such that $z_j^{\alpha_j}=\psi_j^{\alpha_j}(\lambda_2)$,
for $j=1,\ldots,n$.
Then there is a mapping $\tilde\psi\in\Cal O(E,\Bbb C^n)$ such that
$\tilde\psi(E)\subset\subset (\Bbb C_*)^n$, $\psi(\lambda_1)=\tilde\psi(\lambda_1)$, $\tilde\psi(\lambda_2)=z$
and
$$
\gather
\psi_1^{\alpha_1}(\lambda)\cdot\ldots\cdot\psi_n^{\alpha_n}(\lambda)=
\tilde\psi_1^{\alpha_1}(\lambda)\cdot\ldots\cdot\tilde\psi_n^{\alpha_n}
(\lambda),
\quad \lambda\in E,\\
m||\psi_j||_E\leq||\tilde\psi_j||_E\leq M||\psi_j||_E,\quad j=1,\ldots,n
\endgather
$$
where $m,M>0$ depend only on $\delta_1$ and $\alpha$.

\proclaim{Lemma 7} Fix $L^1_1,L^2_1\subset\subset E$,
$L_2\subset\subset \Bbb C_*$ and $\alpha\in(\Bbb R_+)^n$.
Assume that there is $\delta>0$ such that for any $\lambda_1\in L^1_1$,
$\lambda_2\in L^2_1$ we have $m(\lambda_1,\lambda_2)\geq \delta$.

Then there  is $L_2\subset K\subset\subset \Bbb C_*$  such that for any  $z_1,z_2\in L_2$
and for any
$\lambda_1\in L^1_1$, $\lambda_2\in L_1^2$ there is
$\psi\in\Cal O(E,\Bbb C_*)$ with $\psi(\lambda_j)=z_j$, $j=1,2$, and
$\psi(E)\subset K$.

Moreover, there is $\tilde K\subset\subset \Bbb C_*$ such that
for any numbers $z_1,\ldots,z_n\in L_2$,
$w_1,\ldots,w_k\in L_2$, $k<n$ with
$$
|z_1|^{\alpha_1}\cdot\ldots\cdot|z_n|^{\alpha_n}=1
$$
there are functions
$$
\gather
\psi_j\in\Cal O(E,\Bbb C_*),\quad\psi_j(E)\subset \tilde K,\;j=1,\ldots,n,
\quad\psi_1^{\alpha_1}(\lambda)
\cdot\ldots\cdot\psi_n^{\alpha_n}(\lambda)=e^{i\theta},\;
\lambda\in E,\\
\psi_j(\lambda_1)=z_j,\;j=1,\ldots,n,\quad\psi_j(\lambda_2)=w_j,\;
j=1,\ldots,k.
\endgather
$$
\endproclaim
\demo{Proof} For the proof of the first part of the lemma it is sufficient
to prove it for $L^1_1=\{\lambda_1\}$, $L_1^2=\{\lambda_2\}$ with
$m(\lambda_1,\lambda_2)=\delta$. This is so because the general case
one obtains from that special one by composing the functions
with automorphisms of $E$ and dilatation $R\lambda$, where $0\leq R<1$
and as we
see the images of new functions are contained in that of the starting one.

Define
$$
L:=\exp^{-1}(L_2)\cap(\Bbb R\times [0,2\pi))\subset\subset \Bbb C.
$$
Now put
$$
\gather
K:=\{\exp(h(\lambda)):\lambda\in E,\\
\text{ and $h$ is of the type }
h(\lambda)=a\lambda+b, \; a,b\in\Bbb C,
\; h(\lambda_1)=\tilde z_1,\;h(\lambda_2)=\tilde z_2,\;\tilde z_1,\tilde z_2
\in L\}.
\endgather
$$
Observe that $K\subset\subset \Bbb C_*$.
The mappings we are looking for are of the form $\exp\circ h$,
where $h$ is one of the functions appearing in the definition of $K$.

For the proof of the second part of the lemma we put $w_j$ for
$j=k+1,\ldots,n-1$ as any number from $L_2$ and we take mappings
$\psi_1,\ldots,\psi_{n-1}$ as in the first part of the lemma. Define
$$
\psi_n(\lambda):=\frac{e^{i\tilde\theta}}{(\psi_1^{\alpha_1}(\lambda)\cdot\ldots\cdot
\psi_{n-1}^{\alpha_{n-1}}(\lambda))^{1/{\alpha_n}}},
$$
where the branches of powers are chosen arbitrarily and
$\tilde\theta\in\Bbb R$ is chosen so that $\psi_n(\lambda_1)=z_n$.
\qed
\enddemo

\proclaim{Lemma 8} Let $L^1_1,L^2_1,L_2,\delta$ be as in Lemma 7.
Fix $\alpha\in\Bbb
N_*^n$, where $\alpha_j$'s are relatively prime. Then there is
$K\subset\subset \Bbb C_*$ such that  for any mappings $\psi_j\in\Cal O(E,\Bbb
C_*)$, $j=1,\ldots,n$ with
$$
\psi_1^{\alpha_1}\cdot\ldots\cdot\psi_n^{\alpha_n}=1,\quad\lambda\in
E
$$
and $\psi_j(\lambda_1),\psi_j(\lambda_2)\in L_2$, where $\lambda_1\in L_1^1$,
$\lambda_2\in L_1^2$
there are functions $\tilde \psi_j\in \Cal O(E,\Bbb C_*)$ such that
$$
\gather
\tilde\psi_1^{\alpha_1}\cdot\ldots\cdot\tilde\psi_n^{\alpha_n}=1,
\quad\lambda\in E,\\
\tilde\psi_j(\lambda_1)=\psi_j(\lambda_1),\quad
\tilde\psi_j(\lambda_2)=\psi_j(\lambda_2),\quad
\tilde \psi_j(E)\subset K,\quad j=1,\ldots,n.
\endgather
$$
\endproclaim
\demo{Proof}
For     the      proof     put     $z_j:=\psi_j(\lambda_1)$,
$w_j:=\psi_j(\lambda_2)$,   $j=1,\ldots,n$.   From Lemma 7
there   are   $\tilde\psi_j$,
$j=1,\ldots,n-1$ as desired. Put
$$
\tilde\psi_n(\lambda):=\frac{1}{(\tilde\psi_1^{\alpha_1}(\lambda)\cdot
\ldots\cdot\tilde\psi_{n-1}^{\alpha_{n-1}}(\lambda))^{1/\alpha_n}}
$$
We  choose the branch of the power $\frac{1}{\alpha_n}$ so that
$\tilde\psi_n(\lambda_1)=z_n$, note also
that $\tilde\psi_n^{\alpha_n}(\lambda_2)=w_n^{\alpha_n}$, from Remark 6
we may change $\tilde \psi:=(\tilde\psi_1,\ldots,\tilde\psi_n)$
so that all the desired properties are preserved
and, additionally, $\tilde\psi_n(\lambda_2)=w_n$.
\qed
\enddemo
Below we present a lemma, which is a weaker infinitesimal version of
Lemma 7.

\proclaim{Lemma 9} Let $w\in\Bbb C_*$, $X\in\Bbb C$ and $\lambda_1\in E$.
Then there is a mapping $\psi\in\Cal O(E,\Bbb C_*)$ such that
$$
\psi(\lambda_1)=w,\quad \psi^{\prime}(\lambda_1)=X.
$$
Moreover, for given numbers $w_1,\ldots,w_n\in \Bbb C_*$,
$X_1,\ldots, X_k\in\Bbb C$ ($k<n$) and $\alpha\in(\Bbb R_+)^n$,
where $|w_1|^{\alpha_1}\cdot
\ldots\cdot|w_n|^{\alpha_n}=1$ there are mappings
$\psi_j\in\Cal O(E,\Bbb C_*)$, $j=1,\ldots,n$
such that
$$
\gather
\psi_j(\lambda_1)=w_j,\;j=1,\ldots,n,\quad
\psi_j^{\prime}(\lambda_1)=X_j,\; j=1,\ldots,k,\text{ and }\\
\psi_1^{\alpha_1}(\lambda)\cdot\ldots\cdot\psi_n^{\alpha_n}(\lambda)=
e^{i\theta},\quad\lambda\in E.
\endgather
$$
\endproclaim
\demo{Proof}  The first  part goes as in the proof of Lemma 7 (note that
we do not need to specify more, since we do not demand so much about
the mapping $\psi$ as in Lemma 7).
The mapping we are looking for is of the form $\exp(a\lambda+b)$.

For the second part of the  lemma put $X_j$ as any number from $\Bbb
C$  ($j=k+1,\ldots,n-1$). Take  $\psi_j$ as given in
the first  part of  the lemma (for $j=1,\ldots,n-1$)
with $w$ replaced with $w_j$ and $X$ replaced with $X_j$.
Put
$$
\psi_n(\lambda):=\frac{e^{i\tilde\theta}}{(\psi_1^{\alpha_1}(\lambda)\cdot\ldots
\cdot\psi_{n-1}^{\alpha_{n-1}}(\lambda))^{1/\alpha_n}},
$$
where the  branches of powers are  chosen arbitrarily and $\tilde\theta\in\Bbb
R$ is chosen so that $\psi_n(\lambda_1)=w_n$.
\qed
\enddemo
Now we are able to give formulas for the Lempert function and
the Kobayashi--Royden metric for special points.

\proclaim{Lemma 10} Fix  $w\in V_0$.  Let $z\in  D_{\alpha}$ and
$X\in \Bbb C^n$.
Then
$$
\gather
\tilde k_{D_{\alpha}}^*(w,z)=\left(|z_1|^{\alpha_1}\cdot\ldots\cdot
|z_n|^{\alpha_n}\right)^{\frac{1}{\alpha_{j_1}+
\ldots+\alpha_{j_k}}},\\
\kappa_{D_{\alpha}}(w;X)=(|w_1|^{\alpha_1}\cdot\ldots\cdot
|X_{j_1}|^{\alpha_{j_1}}\cdot\ldots\cdot|X_{j_k}|^{\alpha_{j_k}}\cdot
\ldots\cdot |w_n|^{\alpha_n})^{\frac{1}{\alpha_{j_1}+\ldots+\alpha_{j_k}}},
\endgather
$$
where $\Cal J:=\{j\in\{1,\ldots,n\}:w_j=0\}=
\{j_1,\ldots,j_k\}$.
\endproclaim
\demo{Proof}
Without loss of generality we may assume that
$w_1=\ldots=w_k=0$, $w_{k+1},\ldots,w_n\neq 0$,
$n\geq k\geq 1$. We prove both equalities simultanuously.

First we consider the case
$$
z\in \tilde D_{\alpha}
\text{ (respectively, $X_j\neq 0$ for any $j=1,\ldots,k$)}.
$$
Take any $\phi\in\Cal O(\bar E,D_{\alpha})$ such that
$$
\gather
\phi(0)=w,\quad\phi(t)=z
\text{ (respectively, $\phi(0)=w,\quad t\phi^{\prime}(0)=X$),
 for some $t>0$.}
\endgather
$$
We have that
$$
\phi(\lambda)=(\lambda\psi_1(\lambda),\ldots,\lambda\psi_k(\lambda),
\psi_{k+1}(\lambda),\ldots,\psi_n(\lambda)),\quad \psi_j\in
\Cal  O(\bar E,\Bbb C),\;j=1,\ldots,n.
$$
Put
$$
u(\lambda):=\prod_{j=1}^{n}|\psi_j(\lambda)|^{\alpha_j}.
$$
We know that $\log u\in\SH(\bar E)$ and $u\leq 1$ on $\partial E$, so
the maximum principle for subharmonic functions implies that
$u\leq 1$ on $E$. In particular, $u(t)\leq 1$ (respectively, $u(0)\leq 1$),
so
$$
\frac{\prod_{j=1}^{n}|z_j|^{\alpha_j}}{t^{\alpha_1+\ldots+\alpha_k}}\leq 1,
\text{ \big(respectively,
$\frac{\prod_{j=1}^{k}|X_j|^{\alpha_j}\prod_{j=k+1}^{n}|w_j|^{\alpha_j}}
{t^{\alpha_1+\ldots+\alpha_k}}\leq 1$\big),}
$$
which gives us the inequality
$$
t\geq\left(
\prod_{j=1}^{n}|z_j|^{\alpha_j}\right)^{\frac{1}{\alpha_1+\ldots+\alpha_k}},
\text{ \big(respectively, $t\geq\left(
\prod_{j=1}^{k}|X_j|^{\alpha_j}\prod_{j=k+1}^{n}|w_j|^{\alpha_j}\right)
^{\frac{1}{\alpha_1+\ldots+\alpha_k}}$\big)}.
$$

Therefore,
$$
\gather
\tilde k_{D_{\alpha}}^*(w,z)\geq
\left(\prod_{j=1}^{n}|z_j|^{\alpha_j}\right)
^{\frac{1}{\alpha_1+\ldots+\alpha_k}},\\
\text{\big(respectively, $\kappa_{D_{\alpha}}(w;X)\geq
\left(\prod_{j=1}^{k}|X_j|^{\alpha_j}\prod_{j=k+1}^{n}|w_j|^{\alpha_j}\right)
^{\frac{1}{\alpha_1+\ldots+\alpha_k}}$\big)}.
\endgather
$$
To get above the equality put
$$
\gather
t:=\left(\prod_{j=1}^{n}|z_j|^{\alpha_j}\right)^
{\frac{1}{\alpha_1+\ldots+\alpha_k}},\\
\text{\big(respectively,
$t:=\left(\prod_{j=1}^{k}|X_j|^{\alpha_j}\prod_{j=k+1}^{n}|w_j|^{\alpha_j}\right)
^{\frac{1}{\alpha_1+\ldots+\alpha_k}}$\big)}
\endgather
$$
and let us consider the following mapping:
$$
\phi(\lambda):=(\lambda\psi_1(\lambda),\ldots,\lambda\psi_k(\lambda),
\psi_{k+1}(\lambda),\ldots,\psi_n(\lambda)),\quad\lambda\in E,
$$
where $\psi_j\in\Cal O(E,\Bbb C_*)$, $j=1,\ldots,n$,
$\prod_{j=1}^{n}\psi_j(\lambda)^{\alpha_j}=e^{i\theta}$ on $E$
and
$$
\gather
\psi_j(t)=z_j/t,\;j=1,\ldots,k,\quad\psi_j(t)=z_j,\;j=k+1,\ldots,n;\\
\psi_j(0)=w_j,\;j=k+1,\ldots,n, \text{ (see Lemma 7)},\\
\text{\big(respectively, $\psi_j(0)=\frac{X_j}{t}$, $j=1,\ldots,k$,
$\psi_j(0)=w_j$,
$j=k+1,\ldots,n$},\\
\text{$\psi_j^{\prime}(0)=\frac{X_j}{t}$, $j=k+1,\ldots,n$ -- see Lemma 9
\big)}.
\endgather
$$
Then  $\phi\in\Cal O(E,D_{\alpha})$,  $\phi(0)=w$, $\phi(t)=z$
(respectively, $t\phi^{\prime}(0)=X$), which finishes
that case.

We are  remained with the case  $z\in V_0$ (respectively, $X_j=0$
for some $1\leq j\leq k$). If there is
$j$ such that $w_j=z_j=0$ (respectively, $w_j=X_j=0$),
then the mapping
$$
\Bbb C^{n-1}\owns (z_1,\ldots,\check z_j,\ldots,z_n)
\to(z_1,\ldots,0,\ldots,z_n)\in D_{\alpha}
$$
gives us the following
$$
\align
0&=\tilde k^*_{\Bbb C^{n-1}}((w_1,\ldots,\check w_j,\ldots,w_n),
(z_1,\ldots,\check z_j,\ldots,z_n))\geq \tilde k_{D_{\alpha}}^*
(w,z),\\
\text{(}&\text{respectively,}\\
0&=\text{$\kappa_{\Bbb C^{n-1}}((w_1,\ldots,\check w_j,\ldots,w_n);
(X_1,\ldots,\check X_j,\ldots,X_n))\geq \kappa_{D_{\alpha}}(w;X)$)}.
\endalign
$$

Therefore, we are remained only with the Lempert function and then
we may assume that for all $j$ we have $|w_j|+|z_j|>0$.

Define for fixed $\beta>0$ the mapping
$\phi:=(\phi_1,\ldots,\phi_n)$ as follows
$$
\align
\text{if   $w_j=0$},&\text{ then
$\phi_j(\lambda):=\frac{\lambda-\beta}{1-\beta\lambda}
\psi_j(\lambda)$},\\
\text{if $z_j=0$},&\text{ then $\phi_j(\lambda):=
\frac{\lambda+\beta}{1+\beta\lambda}\psi_j(\lambda)$},\\
\text{if $w_jz_j\neq 0$},&\text{ then $\phi_j(\lambda):=\psi_j(\lambda)$},
\endalign
$$
where
$\psi_j\in\Cal O(E,\Bbb C_*)$,
$\prod_{j=1}^{n}\psi_j(\lambda)^{\alpha_j}=e^{i\theta}$ on $E$
and $\phi(\beta)=w$, $\phi(-\beta)=z$ (the values of $\psi_j(\beta)$
and $\psi_j(-\beta)$ are prescribed if only $w_jz_j\neq 0$;
for those $j$ for
which $w_jz_j=0$ only one from the values $\psi_j(\beta)$
and $\psi_j(-\beta)$ is prescribed, more precisely take $j_1$ such that
$z_{j_1}=0$,     then     we     define     $\psi_{j_1}(-\beta)$    so    that
$|\psi_1(-\beta)|^{\alpha_1}\cdot\ldots|\psi_n(-\beta)|^{\alpha_n}=1$;    note
also that
there  is $j_2$  such that  $w_{j_2}=0$, so  $\psi_{j_2}(\beta)$ has  no fixed
value
it is the reason why we are allowed to use Lemma 7). Note also that
$\phi\in \Cal O(E,D_{\alpha})$.
As $\beta>0$ may be chosen arbitrarily small this completes
the proof of the lemma.
\qed
\enddemo
In the next step we will prove a formula for the Lempert function
in the special case
of the domain $D_{(1,\ldots,1)}$. Following (to some extent) the ideas from
\cite{JPZ} and \cite{PZ} we shall propagate the formulas to the general case
using a technic, which could be called a transport of geodesics. Roughly
speaking, the idea relies on transporting the formulas from simpler domains
to more complex ones with the help of 'good' mappings. In \cite{JPZ} and
\cite{PZ} it was the Euclidean ball that was a model domain. In our paper
it is a domain $D_{(1,\ldots,1)}$.

\proclaim{Lemma 11} If $w,z\in V_0$,
then
$$
\tilde k_{D_{(1,\ldots,1)}}^*(w,z)=0.
$$
Assume that $w\in \tilde D_{(1,\ldots,1)}$. Then the following equality holds:
$$
\tilde k_{D_{(1,\ldots,1)}}^*(w,z)=m(w_1\ldots w_n,z_1\ldots z_n)^{1/k},
$$
where
$$
k:=\max\{\#\{j:z_j=0\},1\}.
$$
\endproclaim
\demo{Proof} The first part of the lemma is a consequence of Lemma 10.
Moreover, also the case $z\in V_0$
is a consequence of Lemma 10.

Consider now the case $w,z\in\tilde D_{(1,\ldots,1)}$. We may assume that
$w_1\cdot\ldots\cdot w_n\neq z_1\cdot\ldots\cdot z_n$
(the other case is covered by Remark 4).

Let us consider the following mapping (see Lemma 7):
$$
\phi(\lambda):=\left(\psi_1(\lambda),\ldots,\psi_{n-1}(\lambda),
e^{-i\theta}\lambda\psi_n(\lambda)\right),
$$
where
$$
\gather
\lambda_1:=w_1\cdot\ldots\cdot w_n,\quad\lambda_2:=z_1\cdot\ldots\cdot z_n,\\
\psi_j\in\Cal O(E,\Bbb C_*),\;j=1,\ldots,n,\quad
\psi_1(\lambda)\cdot\ldots\cdot\psi_{n}(\lambda)=e^{i\theta},
\quad \lambda\in E,\\
\psi_j(\lambda_1)=w_j,\quad\psi_j(\lambda_2)=z_j,\quad j=1,\ldots,n-1,
\endgather
$$
(using Lemma 7 we may assume even that $\psi_j(E)\subset K\subset\subset \Bbb C_*$,
$j=1,\ldots,n$ -- compare Remark 12).

Note that
$$
\phi\in\Cal O(E,D_{(1,\ldots,1)}),\quad\phi(\lambda_1)=w,\quad
\phi(\lambda_2)=z.
$$
Therefore, combining these pieces of information with the formula of the
Green function for $D_{(1,\ldots,1)}$ (see Theorem 1) we have:
$$
\multline
m(w_1\cdot\ldots\cdot w_n,z_1\cdot\ldots\cdot z_n)
\geq \tilde k_{D(1,\ldots,1)}^*(w,z)\geq
g_{D_{(1,\ldots,1)}}(w,z)=\\
m(w_1\cdot\ldots\cdot w_n,z_1\cdot\ldots\cdot z_n).
\endmultline
$$
This completes the proof.
\qed
\enddemo

\subheading{Remark 12} From the proof of Lemma 11 we get that for any
$w,z\in\tilde D_{(1,\ldots,1)}$ with $w_1\cdot\ldots\cdot w_n\neq z_1\cdot
\ldots\cdot z_n$ there is a $\tilde k_{D_{(1,\ldots,1)}}$-geodesic
for $(w,z)$, which is of the form
$$
(\psi_1(\lambda),\ldots,\psi_{n-1}(\lambda),
e^{i\theta}\frac{\lambda-\beta}{1-\bar\beta\lambda}\psi_n(\lambda))
$$
with            $\psi_1(\lambda)\cdot\ldots\cdot           \psi_n(\lambda)=1$
and $\psi_j(E)\subset\subset \Bbb C_*$.

\vskip2ex
The domains $D_{\alpha}$ although very regular have not got one property,
which is crucial in the theory of the holomorphically invariant functions;
namely, they are not taut. Therefore, we have no certainty that they admit
$\tilde k_{D_{\alpha}}$-geodesics. However, as Lemma 13 will show, at
least in the rational case and for points, which are 'seperated' by the
Lempert function, it holds. The existence of the geodesics will play
a great  role in  the proof  of the  formula for  the Lempert  function in the
rational case.

\proclaim{Lemma 13} Assume that $\alpha\in\Bbb N_*^n$ and $\alpha_j$'s
are relatively prime.
Let $w,z\in\tilde D_{\alpha}$, $w^{\alpha}\neq z^{\alpha}$.
Then there is a bounded $\tilde k_{D_{\alpha}}$-geodesic $\phi\in\Cal O(E,D_{\alpha})$
for $(w,z)$.
\endproclaim
\demo{Proof} We know that (see Theorem 1)
$$
t:=\tilde                  k_{D_{\alpha}}^*(w,z)\geq
g_{D_{\alpha}}(w,z)=m(w^{\alpha},z^{\alpha})>0;
$$
consequently, there are mappings
$\phi^{(k)}=(\phi^{(k)}_1,\ldots,\phi^{(k)}_n)$, $k=1,2,\ldots$ such that
$$
\phi^{(k)}\in \Cal O(\bar E,D_{\alpha}),\;\phi^{(k)}(0)=w,\;\phi^{(k)}(t_k)=z,
\text{ where $t_k\geq t_{k+1}\to t>0$}.
$$
We have
$$
\phi^{(k)}_j=B_j^{(k)}\psi_j^{(k)},\quad j=1,\ldots,n,
$$
where  $B_j^{(k)}$ is  a Blaschke  product and  $\psi_j^{(k)}\in \Cal O(E,\Bbb
C_*)$.

Put $\psi^{(k)}:=(\psi_j^{(k)})_{j=1}^{n}$.
There  are  two  possibilities  (due   to  maximum  principle  of  subharmonic
functions -- remember about the pseudoconvexity of the domain $D_{\alpha}$):
$$
\gather
\psi^{(k)}(E)\subset D_{\alpha},\tag{2}\\
\psi^{(k)}(E)\subset \partial D_{\alpha}.\tag{3}
\endgather
$$
Below we shall prove that without loss of generality we may reduce our
attention only to the case, which is some kind of generalization of
\thetag{3}.

Take any $k$ such that \thetag{2} is fulfilled. First, notice that the mapping
$\tilde
\psi^{(k)}:=((\psi_1^{(k)})^{\frac{\alpha_1}{\alpha_1\cdot\ldots\cdot\alpha_n}},
\ldots,(\psi_n^{(k)})^{\frac{\alpha_n}{\alpha_1\cdot\ldots\cdot\alpha_n}})
$ is a  mapping from $\Cal  O(E,D_{(1,\ldots,1)})$.
>From Remark 12
there is a $\tilde k_{D_{(1,\ldots,1)}}$-geodesic for
$(\tilde\psi^{(k)}(0),\tilde\psi^{(k)}(t_k))$ of the form
$\mu^{(k)}:=(\hat
\psi_1^{(k)},
\ldots,\hat\psi_{n-1}^{(k)},
e^{i\theta_k}\frac{\lambda-\beta_k}{1-\bar\beta_k\lambda}
\hat\psi_{n}^{(k)})$,
where $\hat\psi_1^{(k)}\cdot\ldots\cdot\hat\psi_n^{(k)}=1$ on $E$ such that
$\mu^{(k)}(0)=\tilde\psi^{(k)}(0)$ and
$\mu^{(k)}(R_kt_k)=\tilde \psi^{(k)}(t_k)$, $\beta_k\in E$, $R_k\leq 1$.

Coming back to the domain $D_{\alpha}$ we see that instead of
considering $\phi^{(k)}$ with the property \thetag{2} we may consider the
mapping
(note  that $\frac{\alpha_1\cdot\ldots\cdot\alpha_n}{\alpha_j}\in\Bbb N$)
$$
\tilde\phi^{(k)}(\lambda):=
\left(B_j^{(k)}(\lambda)(\mu_j^{(k)})
^{\frac{\alpha_1\cdot\ldots\cdot \alpha_n}
{\alpha_j}}(R_k\lambda)\right)_{j=1}^{n},
$$
because $\tilde\phi^{(k)}\in\Cal O(E,D_{\alpha})$,
$\tilde\phi^{(k)}(0)=w$ and $\tilde\phi^{(k)}(t_k)=z$.

Therefore we may assume that (irrespective of which case we start \thetag{2}
or \thetag{3})
$$
\phi_j^{(k)}=B_j^{(k)}\psi_j^{(k)},\quad j=1,\ldots,n,
$$
where $(\psi^{(k)}_1)^{\alpha_1}\cdot\ldots\cdot(\psi^{(k)}_n)^{\alpha_n}=1$ and
$|B_j^{(k)}|\leq 1$ (although we have no longer that $B_j^{(k)}$'s
are the Blaschke products).

Choosing, if necessary, a subsequence we may assume that for all
$j=1,\ldots,n$
$\{B_j^{(k)}\}_{k=1}^{\infty}$ converges
locally uniformly on $E$. Keeping in mind that $\phi^{(k)}(0)=w$ and
$\phi^{(k)}(t_k)=z$ we have in view of Lemma 8
that there is $K\subset\subset \Bbb
C_*$ such that we may assume that
$\psi_j^{(k)}(E)\subset K$ for  any $j,k$
(we may apply Lemma 8 because
$L_2:=\{\psi_j^{(k)}(t_k),\psi_j^{(k)}(0)\}_{j,k}\subset\subset
\Bbb C_*$, which follows from convergence and boundedness of
$\{B_j^{(k)}\}_{k=1}^{\infty}$, the fact that $w_jz_j\neq 0$, $j=1,\ldots,n$
and the equality
$(\psi^{(k)}_1)^{\alpha_1}\cdot\ldots\cdot(\psi^{(k)}_n)^{\alpha_n}=1$),
and then choosing,
if necessary,  a  subsequence  we get that the sequence $\phi^{(k)}$ is
convergent   to  a   mapping  $\phi\in   \Cal  O(E,\bar D_{(\alpha)})$
with $\phi(E)\subset\subset(\Bbb C_*)^n$ such that
$\phi(0)=w$,  $\phi(t)=z$. The  maximum  principle
for subharmonic functions
implies, however, that $\phi(E)\subset D_{\alpha}$. This  completes the proof
of the lemma.
\qed
\enddemo

\subheading{3.  The rational  case --  Proof of  Theorem 2} In the present
section we provide the proof of Theorem 2. Since the theorem
consists of a number of formulas we prove them below one by one. We start
with the Lempert function,
which is basic in the calculation of other functions.

We begin  with a formula  for the  M\"obius  function, which seems  to be very
probable;  nevertheless, we  were not  able  to  find some  references in  the
literature. Its proof  is elementary but it needs  tedious calculations, so we
skip the proof.

\proclaim{Lemma 14} Fix $0<s\leq 1$. Then for any
$\lambda_1\in (0,1)$, $\lambda_2\in E$ we have
$$
m(\lambda_1^{s},\lambda_2^{s})\geq
m(\lambda_1,\lambda_2),
$$
where $\lambda_1^s\in(0,1)$ and the power $\lambda_2^s$ is
chosen  so  that  the left--hand side of the formula is smallest
possible.
\endproclaim

\demo{Proof of formula for $\tilde k_{D_{\alpha}}^*$ in rational case}
The case $w_1\cdot\ldots\cdot w_n=0$ is a consequence of
Lemma 10.
The case $w,z\in\tilde  D_{\alpha}$, $w^{\alpha}=z^{\alpha}$
follows from Remark 4.
We are remained with the case $w,z\in\tilde D_{\alpha}$,
$w^{\alpha}\neq z^{\alpha}$. Due to Lemma 13 there is a bounded
$\tilde k_{D_{\alpha}}$-geodesic $\phi\in\Cal O(E,D_{\alpha})$
for $(w,z)=(\phi(\lambda_1),\phi(\lambda_2))$.
Proceeding
as in the proof of Lemma 13 we may assume that
$$
\phi_j=B_j\psi_j,\quad j=1,\ldots,n,
$$
where $B_j$ is the Blaschke product (up to a constant $|c_j|=1$),
$\psi_j(E)\subset K
\subset\subset \Bbb C_*$ and
$\psi_1^{\alpha_1}\cdot\ldots\cdot\psi_n^{\alpha_n}=1$.
In fact, let us consider
the decomposition of $\phi_j$ as above with the Blaschke product $B_j$. Put
$$
\tilde \psi:=
((\psi_j)^{\frac{\alpha_j}{\alpha_1\cdot\ldots\cdot\alpha_n}})_{j=1}^{n}.
$$
Consider two cases. If $\psi_1^{\alpha_1}\cdot\ldots\cdot\psi_n^{\alpha_n}$
is not constant
on $E$, then $\tilde\psi\in\Cal O(E,D_{(1,\ldots,1)})$ and it is a
$\tilde k_{D_{(1,\ldots,1)}}$-geodesic for $(\tilde\psi(\lambda_1),
\tilde\psi(\lambda_2))$, otherwise, there would be
$\hat\psi\in\Cal O(E,D_{(1,\ldots,1)})$ such that
$\hat\psi(\lambda_1)=\tilde\psi(\lambda_1)$,
$\hat\psi(\lambda_2)=\tilde\psi(\lambda_2)$ and
$\hat\psi(E)\subset\subset D_{(1,\ldots,1)}$, taking then
$\hat\phi(\lambda):=(B_j(\lambda)\hat\psi_j^{\frac{\alpha_1\cdot\ldots\cdot
\alpha_n}{\alpha_j}}(\lambda))_{j=1}^{n}$ we get a mapping
such that $\hat\phi(\lambda_1)=\phi(\lambda_1)$, $\hat\phi(\lambda_2)=
\phi(\lambda_2)$
and $\hat\phi(E)\subset\subset D_{\alpha}$ -- contradiction. Due to
Remark 12 we know that
there is a $\tilde k_{D_{(1,\ldots,1)}}$-geodesic
for $(\tilde \psi(\lambda_1),
\tilde\psi(\lambda_2))=(\mu(\lambda_1),\mu(\lambda_2))$,
where $\hat\psi_1\cdot\ldots\cdot\hat\psi_n=1$ and $\hat\psi_j(E)$'s are
relatively compact in $\Bbb C_*$.
Taking now
$(B_j(\lambda)(\mu_j(\lambda))^{\frac{\alpha_1\cdot\ldots\cdot\alpha_n}
{\alpha_j}})_{j=1}^{n}$
instead    of    $\phi$    we    get    the    desired   property.

In   case $\psi_1^{\alpha_1}\cdot\ldots\cdot\psi_n^{\alpha_n}=e^{i\theta}$,
we     may     assume     that $\psi_j(E)\subset K\subset\subset \Bbb C_*$
for some $K$ because of Lemma 8 (and then without loss of
generality we may assume that
$e^{i\theta}=1$).

Therefore, $\phi(E)$ is contained in some polydisk. Consequently,
$\phi(E)$ is contained in some smooth bounded
pseudoconvex complete Reinhardt
domain $G\subset D_{\alpha}$,
which arises from the domain $D_{\alpha}$ by 'cutting the ends'
and 'smoothing the corners'.
Therefore, $\phi$ is a $\tilde k_G$-geodesic for $(w,z)$.
Using the results of \cite{E}, \cite{Pa} we have that
there are mappings
$h_j\in H^{\infty}(E,\Bbb C)$, $j=1,\ldots,n$ and
$\rho:\partial E\to (0,\infty)$
such that
$$
\frac{1}{\lambda}h_j^*(\lambda)\phi_j^*(\lambda)=\rho(\lambda)
\alpha_j|(\phi^*(\lambda))^{\alpha}|,\quad j=1,\ldots,n,\;
\text{ for almost all $\lambda\in\partial E$}
$$
(we easily exclude the case $(\phi^*(\lambda))^{\alpha}=0$ for $\lambda$ from
some subset of $\partial E$ with non-zero Lebesgue measure).
Using the result of Gentili (see \cite{Ge}) we get that for
some $b_j\in\Bbb C_*$, $j=1,\ldots,n$, $\beta\in E$,
$$
\phi_j(\lambda)h_j(\lambda)=b_j(1-\bar\beta\lambda)(\lambda-\beta),\quad
j=1,\ldots,n,\;\lambda\in E,
$$
where $b_j/\alpha_j=b_k/\alpha_k$, $j,k=1,\ldots,n$.
Consequently, we may take
$$
B_j(\lambda)=c_j
\left(\frac{\lambda-\beta}{1-\bar\beta\lambda}\right)^{r_j},\quad |c_j|=1,
$$
where $r_j\in\{0,1\}$ and not all $r_j$'s are equal to $0$.
Without loss of generality we may assume that $\beta=0$ (we change then only
$\lambda_1$ and $\lambda_2$).

Now we are coming back to the domain $D_{\alpha}$.
We may assume that $r_1=\ldots=r_k=1$ and $r_{k+1}=\ldots=r_n=0$
($1\leq k\leq n$).
We want to have for some $\lambda_1,\lambda_2\in E$ that (without loss of
generality we may assume that $c_j=1$ -- if necessary
we change $w$ and $z$ with the help of rotations of suitable components,
so the Lempert function does not change)
$$
\gather
\lambda_1\psi_j(\lambda_1)=w_j,\; j=1,\ldots,k,\quad
\psi_j(\lambda_1)=w_j,\;j=k+1,\ldots,n,\\
\lambda_2\psi_j(\lambda_2)=z_j,\; j=1,\ldots,k,\quad
\psi_j(\lambda_2)=z_j,\;j=k+1,\ldots,n.
\endgather
$$
Taking the $\alpha_j$-th power and multiplying the equalities we get
that
$$
\lambda_1^{\alpha_1+\ldots+\alpha_k}=w^{\alpha},\quad
\lambda_2^{\alpha_1+\ldots+\alpha_k}=z^{\alpha}.
$$
The formulas above describe all possibilities, which
may deliver us the candidates for the realization of the
Lempert function.
Now for all possible pairs of numbers
$\lambda_1,\lambda_2$ given as above we find mappings, which map
$\lambda_1$ and $\lambda_2$ in $w$ and $z$.
Note that there are mappings $\psi_j\in\Cal O(E,\Bbb C_*)$,
$j=2,\ldots,n$ such that (see Lemma 7)
$$
\gather
\psi_j(\lambda_1)=
\frac{w_j}{(w^{\alpha})^{\frac{1}{\alpha_1+\ldots+\alpha_k}}}=
\frac{w_j}{\lambda_1},
\quad j=2,\ldots,k,\\
\psi_j(\lambda_2)=\frac{z_j}
{(z^{\alpha})^{\frac{1}{\alpha_1+\ldots+\alpha_k}}}=\frac{z_j}{\lambda_2},
\quad j=2,\ldots,k,\\
\psi_j(\lambda_1)=w_j,\quad j=k+1,\ldots,n,\\
\psi_j(\lambda_2)=z_j,\quad j=k+1,\ldots,n.
\endgather
$$
Define also
$$
\psi_1(\lambda):=
\frac{1}{(\psi_2^{\alpha_2}(\lambda)\cdot\ldots\cdot
\psi_{n}^{\alpha_{n}}(\lambda))^
{\frac{1}{\alpha_1}}},\quad\lambda\in E.
$$
Put
$$
\phi(\lambda):=(\lambda\psi_1(\lambda),\ldots,\lambda\psi_k(\lambda),
\psi_{k+1}(\lambda),\ldots,\psi_n(\lambda)).
$$
The $\frac{1}{\alpha_1}$-st root in definition of
$\psi_1$  is chosen so that
$\phi_1(\lambda_1)=w_1$,
and  we  know  that  $\phi_1^{\alpha_1}(\lambda_2)=z_1^{\alpha_1}$.
One may also easily
verify that $\phi(\lambda_1)=w$ and $\phi_j(\lambda_2)=z_j$, $j=2,\ldots,n$,
which, however, in view of Lemma 5 shows that there is also a mapping
$\tilde\phi\in \Cal O(E,D_{\alpha})$
such that $\tilde\phi(\lambda_1)=w$, $\tilde\phi(\lambda_2)=z$.
Therefore we have proved that
$$
\tilde
k_{D_{\alpha}}^*(w,z)=\min\{m(\lambda_1,\lambda_2):\lambda_1,
\lambda_2\in E,
\lambda_1^{\alpha_{j_1}+\ldots+\alpha_{j_k}}=w^{\alpha},
\lambda_2^{\alpha_{j_1}+\ldots+\alpha_{j_k}}=z^{\alpha}\},
$$
where    the     minimum    is    taken    over     all    possible    subsets
$\{j_1,\ldots,j_k\}\subset\{1,\ldots,n\}$.
And now Lemma 14 finishes the proof (remark that without loss of generality we
may assume that $w_j>0$, $j=1,\ldots,n$).
\qed
\enddemo
\demo{Proof of the formula for $k_{D_{\alpha}}^*$ in rational case}
Note that $\tanh^{-1}$ of
the desired formula is equal to $\tanh^{-1}$
of the Lempert function off the axis, satisfies the
triangle  inequality  and  is  continuous.
The  definition  of  the  Kobayashi
pseudodistance and its continuity (see \cite{JP 2}) finish the proof.
\qed
\enddemo

To finish the proof we are remained only with the problem of computing
the Kobayashi--Royden pseudometric $\kappa_{D_{\alpha}}$. We get that formula
from that of the Kobayashi pseudodistance. But to see that we have to
define an operator, which connects these both functions.

Following M. Jarnicki and  P. Pflug (see \cite{JP 2}) for a  domain
$D\subset \Bbb C^n$ we define the following function
$$
\frak D k_D(w;X):=
\limsup_{\lambda\not\to 0}\frac{k_{D_{\alpha}}^*(w,w+\lambda X)}
{|\lambda|},\quad w\in D, X\in\Bbb C^n.
$$
The function defined above differs from that in \cite{JP 2},
nevertheless,
since our version is not larger than that from \cite{JP 2}
the inequality below, which is crucial for our considerations,
remains true
$$
\frak D k_D(w;X)\leq\kappa_D(w;X),\quad w\in D,X\in \Bbb C^n.\tag{4}
$$

\proclaim{Lemma 15}  Let  $\alpha\in\Bbb   N_*^n$,  where  $\alpha_j$'s  are
relatively prime. Then
$$
\multline
\frak D k_{D_{\alpha}}(w;X)=\\
\gamma_E\left(\left(\prod_{j=1}^n|w_j|^{\alpha_j}\right)
^{\frac{1}{\min\{\alpha_k\}}},
\left(\prod_{j=1}^n|w_j|^{\alpha_j}\right)
^{\frac{1}{\min\{\alpha_k\}}}\frac{1}{\min\{\alpha_k\}}
\sum_{j=1}^n\frac{\alpha_j X_j}
{w_j}\right),
\endmultline
$$
$w\in\tilde  D_{\alpha},X\in  \Bbb
C^n$.
\endproclaim
\demo{Proof}  Without   loss  of  generality  we   may  assume  that  $w_j>0$,
$j=1,\ldots,n$ and
$\alpha_n=\min\{\alpha_k\}$.
Using the formula for $k_{D_{\alpha}}^*$ we get
$$
\frak D k_{D_{\alpha}}(w;X)=
\limsup_{\lambda\not\to 0}\frac{\left|\prod_{j=1}^n(w_j+\lambda X_j)
^{\alpha_j/\alpha_n}-
\prod_{j=1}^nw_j^{\alpha_j/\alpha_n}\right|}
{\left|1-\prod_{j=1}^n(w_j+\lambda X_j)^{\alpha_j/\alpha_n}
\prod w_j^{\alpha_j/\alpha_n}\right||\lambda|}.\tag{5}
$$
Applying the Taylor formula we get for $\lambda$ close to $0$
$$
(w_j+\lambda X_j)^{\alpha_j/\alpha_n}=w_j^{\alpha_j/\alpha_n}+
\frac{\alpha_j}{\alpha_n}w_j^{\alpha_j/\alpha_n}\frac{\lambda
X_j}{w_j}+\varepsilon_j(\lambda), \quad j=1,\ldots,n,
$$
where $\frac{\varepsilon_j(\lambda)}{\lambda}
\to 0$ as $\lambda\to 0$. Substituting the last equalities in \thetag{5}
we get that
$$
\frak D k_{D_{\alpha}}(w;X)=\limsup_{\lambda\not\to 0}
\frac{\left(\prod_{j=1}^n|w_j^{\alpha_j}|^{1/\alpha_n}\right)
|\lambda|\left|\sum_{j=1}^{n}\frac{\alpha_jX_j}{\alpha_nw_j}
\right|}
{\left(1-\prod_{j=1}^n|w_j|^{2\alpha_j/\alpha_n}\right)|\lambda|},
$$
which equals the desired value.
\qed
\enddemo

\demo{Proof of the formula for $\kappa_{D_{\alpha}}$  in rational case}
If $\Cal J\neq \emptyset$, then, in view of Lemma 10 we are done.
The case $\sum_{j=1}^{n}\frac{\alpha_jX_j}{w_j}=0$ follows from Remark 4.

Take $w\in\tilde D_{\alpha}$. Without loss of generality we may assume that
$w_j\in \Bbb R_+$, $j=1,\ldots,n$ and $\alpha_n=\min\{\alpha_j\}$.
Below, for $X\in \Bbb C^n$,
$\sum_{j=1}^{n}\frac{\alpha_jX_j}{w_j}\neq 0$ we shall construct a mapping
$\phi\in\Cal O(E,D_{\alpha})$ such that
$$
\phi(\lambda_1)=w,\quad t\phi^{\prime}(\lambda_1)=X,
$$
where
$\lambda_1:=(w_1^{\alpha_1}\cdot\ldots\cdot w_n
^{\alpha_n})^{1/\alpha_n}>0$, $t:=
(w_1^{\alpha_1}\cdot\ldots\cdot w_n
^{\alpha_n})^{1/\alpha_n}\sum_{j=1}^{n}\frac{\alpha_jX_j}
{\alpha_n w_j}$.

Note that the existence of such a $\phi$ would finish the proof
because  of Lemma 15 and \thetag{4}.

Define the mapping
$$
\phi(\lambda):=(\psi_1(\lambda),\ldots,\psi_{n-1}(\lambda),
\frac{\lambda}
{(\psi_1^{\alpha_1}(\lambda)\cdot\ldots\cdot\psi_{n-1}^{\alpha_{n-1}}(\lambda))^
{1/\alpha_n}}),
$$
where (see Lemma 9)
$$
\psi_j(\lambda_1)=w_j,\;j=1,\ldots,n-1,
\quad t\psi_j^{\prime}(\lambda_1)=X_j,\;
j=1,\ldots,n-1.
$$
We choose the $\frac{1}{\alpha_n}$-th power so that $\phi_n(\lambda_1)=
w_n$, after some elementary transformation we get that
$$
t\phi_n^{\prime}(\lambda_1)=X_n,
$$
which finishes the proof.
\qed
\enddemo

\subheading{4. The irrational case -- Proof of Theorem 3} As in rational case
we start with the proof of the formula of the Lempert function. First, we make
use of the special properties of the domains of irrational type to get:

\proclaim{Lemma 16} Let $\alpha$ be of irrational type. Then for
any $w,z\in D_{\alpha}$
$$
\tilde k_{D_{\alpha}}^*(w,z)=\tilde k_{D_{\alpha}}^*(\tilde w,\tilde z),
\quad
\tilde w\in T_w,\; \tilde z\in T_z.
$$
\endproclaim
\demo{Proof} Certainly it is enough to prove that
$$
\tilde k_{D_{\alpha}}^*(w,z)=\tilde k_{D_{\alpha}}^*(w,\tilde z),
\text{ whenever $\tilde z\in T_z$}.
$$
Assume that
$$
\tilde k_{D_{\alpha}}^*(w,\tilde z_1)<
\tilde k_{D_{\alpha}}^*(w,\tilde z_2)=:\varepsilon\tag{6}
$$
for some $\tilde z_1,\tilde z_2\in T_z$.
Then in view of Lemma 5
$$
\tilde k_{D_{\alpha}}^*(w,\tilde z)=\varepsilon\tag{7}
$$
for all $\tilde z\in T_{\tilde z_2,\alpha}$. Because of \thetag{1} we have
that $\tilde z_1\in T_z=T_{\tilde z_2}=\bar  T_{\tilde z_2,\alpha}$.
The  last statement contradicts,  in
connection with \thetag{6} and \thetag{7},
however, the upper-semicontinuity of the Lempert function.
\qed
\enddemo
As an immediate corollary of Lemma 16 we get

\proclaim{Corollary 17} Let $\alpha$ be of irrational type, then for any
$z\in D_{\alpha}$
$$
\tilde k_{D_{\alpha}}^*(z,\tilde z)=0\text{ for any $\tilde z\in T_z$.}
$$
\endproclaim

\demo{Proof of the formula for $\tilde k_{D_{\alpha}}^*$ in irrational case}
The case $\Cal J\neq \emptyset$ is covered by Lemma 10.
Consider now the remaning case. In view of Lemma 16 we have that
$$
\tilde k_{D_{\alpha}}^*(w,z)=\tilde k_{D_{\alpha}}^*((|w_1|,\ldots,|w_n|),
(|z_1|,\ldots,|z_n|)).
$$

%Therefore, for the remaining part of the proof we assume that
%$w,z\in\Bbb R_+^n\cap D_{\alpha}$.

Let us choose a sequence $\{\alpha^{(k)}\}_{k=1}^{\infty}
\subset(\Bbb Q_+)^n$ such that
$$
\alpha^{(k)}\to\alpha.
$$
First notice  that in view  of Theorem 2  we know that
if $x,y\in(\Bbb R_+)^n\cap D_{\alpha^{(k)}}$, then
$$
\tilde k_{D_{\alpha^{(k)}}}^*(x,y)=
m((x_1^{\alpha_1^{(k)}}\cdot\ldots\cdot x_n^{\alpha_n^{(k)}})
^{\frac{1}{\min\{\alpha_j^{(k)}\}}},
(y_1^{\alpha_1^{(k)}}\cdot\ldots\cdot y_n^{\alpha_n^{(k)}})
^{\frac{1}{\min\{\alpha_j^{(k)}\}}}).\tag{8}
$$
We may assume that $\min\{\alpha_j\}=\alpha_n$ and $\min\{\alpha^{(k)}_j\}=
\alpha^{(k)}_n$.
First we prove that
$$
\tilde k_{D_{\alpha}}^*(w,z)\geq
m((|w_1|^{\alpha_1}\cdot\ldots\cdot|w_n|^{\alpha_n})
^{1/\alpha_n},
(|z_1|^{\alpha_1}\cdot\ldots\cdot|z_n|^{\alpha_n})
^{1/\alpha_n}).
$$
Suppose it does not hold, so there is a mapping
$\phi\in \Cal O(\bar E,D_{\alpha})$ such that $\phi(\lambda_1)=
(|w_1|,\ldots,|w_n|)$,
$\phi(\lambda_2)=(|z_1|,\ldots,|z_n|)$ and
$$
m(\lambda_1,\lambda_2)<m((|w_1|^{\alpha_1}\cdot\ldots\cdot|w_n|^{\alpha_n})
^{1/\alpha_n},
(|z_1|^{\alpha_1}\cdot\ldots\cdot|z_n|^{\alpha_n})
^{1/\alpha_n}).
$$
Then we may choose $k$ so large that $\phi(E)\subset D_{\alpha^{(k)}}$ and
\linebreak
$m(\lambda_1,\lambda_2)<m((|w_1|^{\alpha^{(k)}_1}\cdot\ldots\cdot|w_n|
^{\alpha^{(k)}_n})^{1/\alpha^{(k)}_n},
(|z_1|^{\alpha^{(k)}_1}\cdot\ldots\cdot|z_n|^{\alpha^{(k)}_n})
^{1/\alpha^{(k)}_n})$,
which, however, contradicts \thetag{8}.

To get the equality consider the
mapping $\phi(\lambda):=(\psi_1(\lambda),\ldots,\psi_{n-1}(\lambda),
\lambda\psi_n(\lambda))$,
where (see Lemma 7)
$$
\gather
\psi_j\in\Cal O(E,\Bbb C_*),\quad j=1,\ldots,n-1,\\
\lambda_1:=(|w_1|^{\alpha_1}\cdot\ldots\cdot |w_n|^{\alpha_n})
^{\frac{1}{\alpha_n}}>0,
\lambda_2:=(|z_1|^{\alpha_1}\cdot\ldots\cdot|z_n|^{\alpha_n})
^{\frac{1}{\alpha_n}}>0;\\
\psi_j(\lambda_1)=|w_j|,
\quad\psi_j(\lambda_2)=|z_j|,\;j=1,\ldots,n-1,\\
\endgather
$$
Define also
$$
\psi_n(\lambda):=
\frac{1}{(\psi_1^{\alpha_1}(\lambda)\cdot\ldots\cdot
\psi_{n-1}^{\alpha_{n-1}}(\lambda))^
{\frac{1}{\alpha_n}}},\quad\lambda\in E.
$$
The $\frac{1}{\alpha_n}$-th root is chosen so that
$\phi_n(\lambda_1)=|w_n|$. One may also easily check from
the form of $\psi_j$'s in the proof of Lemma 7
that then $\phi_n(\lambda_2)>0$, so
$\phi_n(\lambda_2)=|z_n|$. This completes the proof.
\qed
\enddemo
Identically as in the rational case we have:
\demo{Proof of the formula for $k_{D_{\alpha}}^*$ in irrational case}
Note that $\tanh^{-1}$ of the desired formula satisfies the
triangle inequality and coincides with the $\tanh^{-1}$ of
Lempert function off the axis.
The  continuity  of  the  Kobayashi  pseudodistance  (see \cite{JP 2})
as well as the definition of the Kobayashi pseudodistance finish the proof.
\qed
\enddemo
Having the formula for the Lempert function we get

\demo{Proof of the formula for $g_{D_{\alpha}}$ in the irrational case}

{\bf Case I.} $\Cal J=\emptyset$.

Corollary 16 implies that
$$
g_{D_{\alpha}}(w,z)=0\text{ for any $z\in T_w$}.
$$
Maximum    principle    for    plurisubharmonic    functions    (applied    to
$g_{D_{\alpha}}(w,\cdot)$) implies that
$$
g_{D_{\alpha}}(w,z)=0\text{ for any $z$ with $|z_j|\leq |w_j|$},
$$
which, however,  means that $g_{D_{\alpha}}(w,\cdot)$  vanishes on a  set with
non-empty  interior  (remember  that   $w_1\cdot\ldots\cdot  w_n\neq  0$)  but
$g_{D_{\alpha}}(w,\cdot)$  is  logarithmically  plurisubharmonic,  so  it must
vanish on $D_{\alpha}$.

{\bf Case II.} $\Cal J\neq \emptyset$.

This case is a simple consequence of Lemma 10,
the inequality
$g\leq \tilde k^*$, definition of the Green function and the fact that
the function $(|z_1|^{\alpha_1}\cdot\ldots\cdot|z_n|^{\alpha_n})
^{\frac{1}{\alpha_{j_1}+\ldots+\alpha_{j_k}}}$
is logarithmically plurisubharmonic on $D_{\alpha}$.
\qed
\enddemo
\demo{Proof of the formula for $A_{D_{\alpha}}$ in irrational case}
The result follows  from the formula for the Green  function
and definition of the Azukawa pseudometric.
\qed
\enddemo
And now similarly as in the rational case we finish up the proof by showing
the formula for $\kappa_{D_{\alpha}}$.

\proclaim{Lemma 18} Let $\alpha$ be of irrational type. Then
$$
\multline
\frak D k_{D_{\alpha}}(w;X)=\\
\gamma_E\left(\left(\prod_{j=1}^n|w_j|^{\alpha_j}\right)
^{\frac{1}{\min\{\alpha_k\}}},
\left(\prod_{j=1}^n|w_j|^{\alpha_j}\right)
^{\frac{1}{\min\{\alpha_k\}}}\frac{1}{\min\{\alpha_k\}}
\sum_{j=1}^n\frac{\alpha_j X_j}
{w_j}\right),
\endmultline
$$
for $w\in\tilde  D_{\alpha},X\in  \Bbb C^n$.
\endproclaim
\demo{Proof} Without loss of generality we may assume that
$\alpha_n=\min\{\alpha_k\}$.
The formula for the Kobayashi pseudodistance gives us
$$
\frak D k_{D_{\alpha}}(w;X)=
\limsup_{\lambda\not\to 0}\frac{\left|\prod_{j=1}^n|w_j+\lambda X_j|
^{\alpha_j/\alpha_n}-
\prod_{j=1}^n|w_j|^{\alpha_j/\alpha_n}\right|}
{\left|1-\prod_{j=1}^n|w_j+\lambda X_j|^{\alpha_j/\alpha_n}
\prod|w_j|^{\alpha_j/\alpha_n}\right||\lambda|}.\tag{9}
$$
Note that $\alpha_j/\alpha_n\geq 1$. Therefore,
applying the Taylor formula we get, for $\lambda$ close to $0$,
$$
|w_j+\lambda X_j|^{\alpha_j/\alpha_n}=|w_j|^{\alpha_j/\alpha_n}+
\frac{\alpha_j}{\alpha_n}|w_j|^{\alpha_j/\alpha_n}
\left|\Re\left(\frac{\lambda
X_j}{w_j}\right)\right|+\varepsilon_j(\lambda), \quad j=1,\ldots,n,
$$
where $\frac{\varepsilon_j}{\lambda}\to 0$
as $\lambda\to 0$. Substituting the last equalities in \thetag{9}
we get that
$$
\frak D k_{D_{\alpha}}(w;X)=\limsup_{\lambda\not\to 0}
\frac{\prod_{j=1}^n(|w_j|^{\alpha_j})^{1/\alpha_n}\Re
\left(\lambda\left(\sum_{j=1}^{n}\frac{\alpha_jX_j}{\alpha_nw_j}
\right)\right)}{\left(1-\prod_{j=1}^n|w_j|^{2\alpha_j/\alpha_n}\right)
|\lambda|},
$$
which equals the desired value.
\qed
\enddemo

\demo{Proof of the formula for $\kappa_{D_{\alpha}}$ in irrational case}
If $\Cal J\neq \emptyset$, then, in view of Lemma 10, we are done.
Also the case $\sum_{j=1}^{n}\frac{\alpha_jX_j}{w_j}=0$ follows
from Remark 4.
Below we deal with the remaining cases.

Take $w\in\tilde D_{\alpha}$.
Without loss of generality we may assume that
$w_j\Bbb\in \Bbb R_+$, $j=1,\ldots,n$ and $\alpha_n=\min\{\alpha_j\}$.
Below, for $X\in \Bbb C^n$ with
$\sum_{j=1}^{n}\frac{\alpha_jX_j}{w_j}\neq 0$ we shall construct a mapping
$\phi\in\Cal O(E,D_{\alpha})$ such that
$$
\phi(\lambda_1)=w,\quad t\phi^{\prime}(\lambda_1)=X,
$$
where
$\lambda_1:=(w_1^{\alpha_1}\cdot\ldots\cdot w_n
^{\alpha_n})^{1/\alpha_n}>0$, $t:=
(w_1^{\alpha_1}\cdot\ldots\cdot w_n
^{\alpha_n})^{1/\alpha_n}\sum_{j=1}^{n}\frac{\alpha_jX_j}
{\alpha_n w_j}$.

Note that the existence of such a $\phi$ would finish the proof
because  of Lemma 18 and \thetag{4}.

Define the mapping
$$
\phi(\lambda):=\left(\psi_1(\lambda),\ldots,\psi_{n-1}(\lambda),
\frac{\lambda}
{(\psi_1^{\alpha_1}(\lambda)\cdot\ldots\cdot\psi_{n-1}^{\alpha_{n-1}}(\lambda))^
{1/\alpha_n}}\right),
$$
where (see Lemma 9)
$$
\psi_j(\lambda_1)=w_j,\;j=1,\ldots,n-1,
\quad t\psi_j^{\prime}(\lambda_1)=X_j,\;
j=1,\ldots,n-1.
$$
We choose the $1/\alpha_n$-th power so that
$\phi_n(\lambda_1)=w_n$.
After some elementary transformation we get that
$$
t\phi_n^{\prime}(\lambda_1)=X_n,
$$
which finishes the proof.
\qed
\enddemo

\subheading{5. Some applications}
Having proven the formulas for the invariant
functions  for  the  elementary  Reinhardt
domains  we  may  formulate  some
conclusions, which follow  from them.
They show how  irregularly the functions
can behave although the domains considered are very regular.

For a given domain $D\subset \Bbb C^n$ we define a relation $\Cal R$
on  $D$ as  follows: $w\Cal  R z$  for $w,z\in D$ if $k^*_D(w,z)=0$. In
\cite{Ko2}, S. Kobayashi asked a question whether the quotient $D/\Cal R$
has always
a complex structure. The answer is 'no',
however the examples showing this
are artificial  (see \cite{Ko 1}, page 130 also \cite{HD} and \cite{Gi}).
>From  Theorem 3 we know  that if
$\alpha$  is of  irrational type,
then $D_{\alpha}/\Cal  R$ is equal to $[0,1)$. This
gives the first very simple example of a very regular domain, for
which the answer to the above question is 'no'.

One may consider some generalizations of the
Carath\'eodory pseudodistance, called the $k$-th M\"obius function,
denoted by $m^{k}$ (for $k=1,2,\ldots$)
%Analoguously one may define a generalization
%of the Carath\'eodory--Reiffen pseudometric, so called
%$k$-th Carath\'eodory--Reiffen pseudometric, denoted by $\gamma^k$
(for definitions see \cite{JP 2}).
S. Nivoche (see \cite{N}) has proved that if a domain is strictly
hyperconvex, then the functions $m^k$ tend to
$g$. One may easily verify that if $\alpha$ is of
irrational type, then all the $m^k$'s vanish on
$D_{\alpha}\times D_{\alpha}$. Therefore
we see that no such convergence takes place in domains $D_{\alpha}$
($\alpha$ of irrational type), so one may
not expect a similar result as in  \cite{N}
in the class of Reinhardt complete
pseudoconvex domains.

In  general, the Lempert  function seems to  be very distant
from the
Green  function. The definition of the Kobayashi pseudodistance makes
an impression that the Kobayashi pseudodistance should  be larger or equal
than the Green  function.
Nevertheless, if $\alpha\in\Bbb N^n_*$  is such that
all $\alpha_j$'s  are relatively prime and  $\min\{\alpha_j\}=1$,
then we have the following inequalities (see Theorem 2):
$$
c_{D_{\alpha}}^*\equiv    k_{D_{\alpha}}^*\leq    g_{D_{\alpha}}\leq    \tilde
k_{D_{\alpha}}^*,\quad k_{D_{\alpha}}^*\neq g_{D_{\alpha}}.
$$

In  the  papers  \cite{Pa}  and  \cite{L1}  a  notion  of stationary maps was
introduced and studied.
In the  class of strongly convex domains these mappings  are exactly the
$\tilde k$-geodesics.
In case of  strongly pseudoconvex domains godesics must
be  necessarily  stationary  maps.
One  knows  that  in  general  the inverse
implication  does not  hold (see  \cite{Pa} and  \cite{PZ}).
>From  the proof of Theorem 2
we may construct also other examples  disproving the implication. One may
find even domains, which are strongly pseudoconvex
(one produces  them by cutting the 'ends' and then smoothing the corners).

\subheading{6. Open  problems} It  would be  interesting to  find formulas of all
invariant functions discussed above for  domains of the following more general
type:
$$
D_{\alpha^1}\cap\ldots\cap D_{\alpha^k}\cap ((R_1 E)\times\ldots\times(R_n E)),
$$
where $\alpha^j\in(\Bbb R_+)^n$, $j=1,\ldots,k$.

\enddocument